\documentclass[a4paper]{jpconf}
\usepackage{graphicx}
\usepackage{floatrow}

\begin{document}
\title{The Unexpected Fractal Signatures in Fibonacci chains}

\author{Fang Fang, Raymond Aschheim and Klee Irwin}

\address{Quantum Gravity Research, Los Angeles, CA, U.S.A.}

\ead{Fang@QuantumGravityResearch.org}

\begin{abstract}

Quasicrystals are fractal due to their self similar property. In this paper, a new cycloidal fractal signature possessing the cardioid shape in the Mandelbrot set is presented in the Fourier space of a Fibonacci chain with two lengths, $L$ and $S$, where $L/S=\phi$. The corresponding pointwise dimension is $0.7$. Various modifications, such as truncation from the head or tail, scrambling the orders of the sequence, and changing the ratio of the $L$ and $S$, are done on the Fibonacci chain. The resulting patterns in the Fourier space show that that the fractal signature is very sensitive to changes in the Fibonacci order but not to the $L/S$ ratio.
\end{abstract}

\section{Introduction}
Quasicrystals possess exotic and sometimes anomalous properties that have interested the scientific community since their discovery by Shechtman in 1982 \cite{shechtman1984metallic}. Of particular interest in this manuscript is the self similar property of quasicrystals that makes them fractal. Historically, research on the fractal aspect of quasicrystalline properties has revolved around spectral and wave function analysis \cite{GrimmQCHamiltonian, RepetowiczPTHamiltonians, MaceFiChainFractal2016}. Mathematical investigation \cite{Ramachandrarao2000, Yudin2001} of the geometric structure of quasicrystals is less represented in the literature than experimental work.

In this paper, a new framework for analyzing the fractal nature of quasicrystals is introduced. Specifically, the fractal properties of a one-dimensional Fibonacci chain and its derivatives are studied in the complex Fourier space. The results may also be found in two and three dimensional quasicrystals that can be constructed using a network of one dimensional Fibonacci chains \cite{Fang2015}.

\section{The fractal signature of the Fibonacci chain in Fourier space}

The Fibonacci chain is a quasiperiodic sequence of short and long segments where the ratio between the long and the short segment is the golden ratio  \cite{levine1986quasicrystals}. It is an important 1D quasicrystal that uniquely removes the arbitrary closeness in the non-quasicrystalline grid space of a quasicrystal. And as a result the non-quasicrystalline grid space is converted into a quasicrystal as well \cite{Fang2015}. The Fourier representation of a Fibonacci chain is given as follows:

 \begin{equation} \label{eq:1}
z_s=\frac{1}{\sqrt{n}}\sum _{r=1}^n u_re^{\frac{2 \pi  i (r-1) (s-1)}{n}}, s=1,2,...n,
\end{equation}
where 
 \begin{equation} \label{eq:2}
u_r=(\phi -1) (\lfloor (r+1) \phi \rfloor -\lfloor r \phi \rfloor -2)+\phi =(\phi -1) (\lfloor (r+1) \phi \rfloor -\lfloor r \phi \rfloor )+(2-\phi ), r=1,2,...n
\end{equation}
is a Fibonacci chain of length $n$ with units of length $\phi$ and $1$. Here $\phi= \frac{ \sqrt{5} +1 }{2} $ is the golden ratio.

Note that the Fibonacci chain can also be generated using substitution rules \cite{levine1986quasicrystals} where $m$ iterations of the substitution method produces a Fibonacci chain of length equal to the $m$th Fibonacci number.

\begin{figure}[h!]
  \centering
  \includegraphics[width=0.8\textwidth,clip]{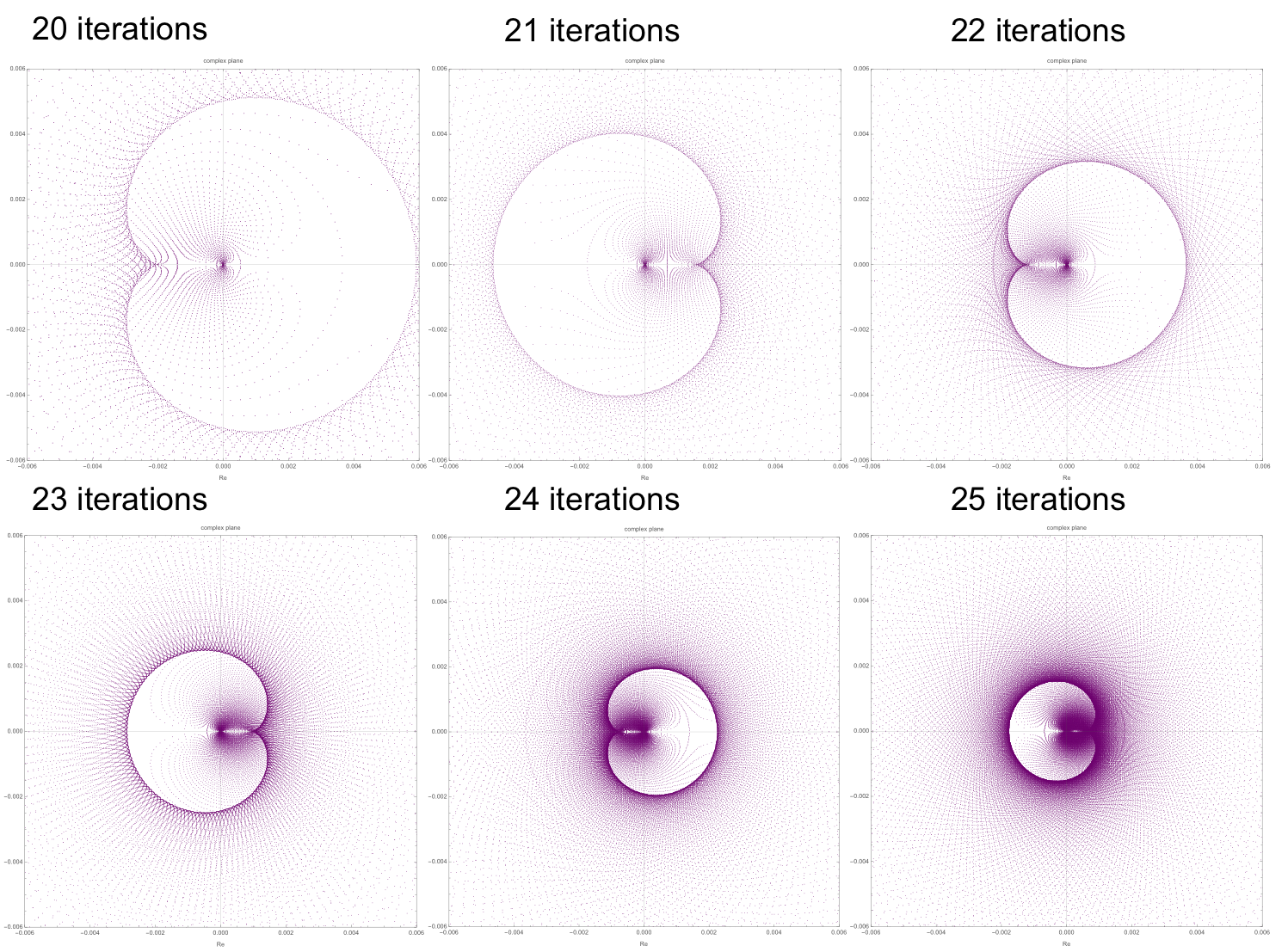}
  \caption{The Fourier space representation of a Fibonacci chain is shown here using the substitution method with different iterations, starting with $20$. The horizontal axis represents the real part and the vertical axis represents the imaginary part of the Fourier coefficients.}
  \label{f1}
\end{figure}

Consecutive iterations of the Fibonacci chain are plotted in the Fourier space in Fig. \ref{f1}. The main cardioid (heart shape) in the Mandelbrot set \cite{douady1984etude} appeared at the center of each plot. It is known that the denominators of the periods of the circular bulbs at sequential scales in the Mandelbrot set follow the Fibonacci number sequence \cite{DevaneyBU}. The plots reveal that the size of the cardioid decreases with successive iterations, and the scaling factor for the size change is approximately $1/\sqrt{\phi}$. Moreover, the orientation of the cardioid undergoes mirror flips with every odd and even numbered iteration. Upon magnification, the fractal nature becomes apparent near the real axis. The fractal dimension\footnote{pointwise dimension \cite{Strogatz1994} } of the complex Fourier representation of the Fibonacci chain is approximately $0.7$. 

\begin{figure}[h!]
  \centering
  \includegraphics[width=0.8\textwidth,clip]{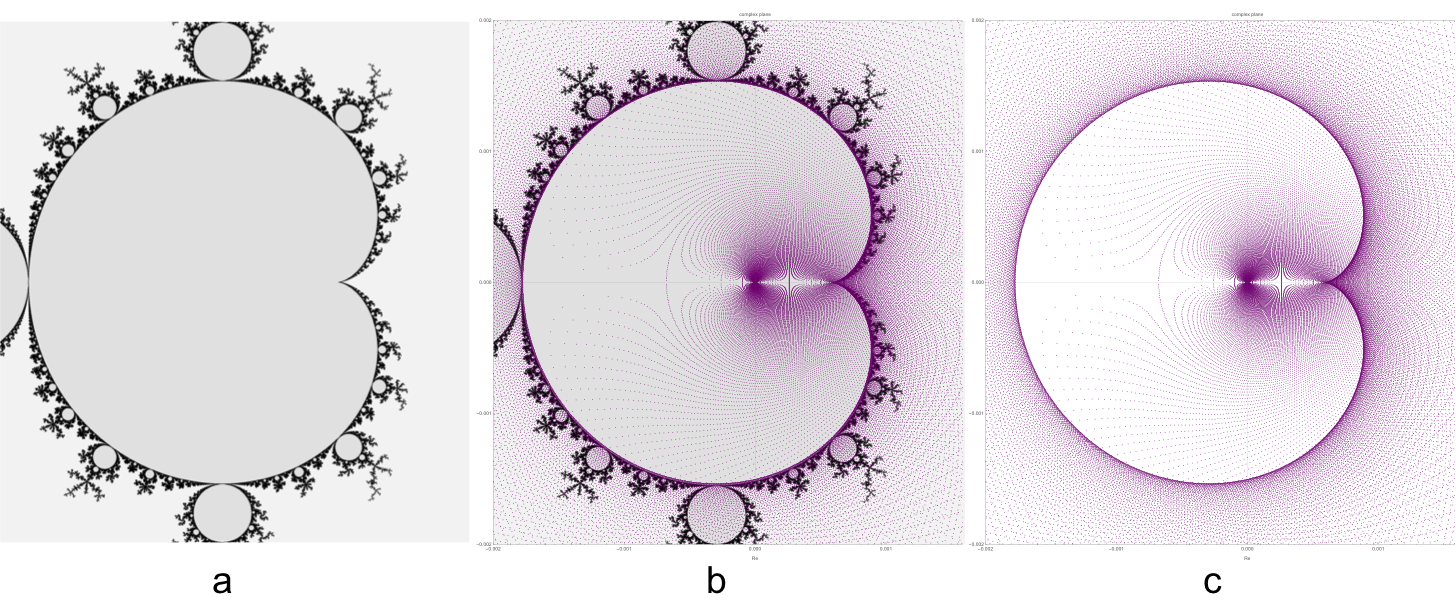}
  \caption{a) The Mandelbrot set, c) the fractal structure in the Fourier space of a Fibonacci chain with 25 iterations and b) overlays \textit{a} and \textit{c} together.}
  \label{f1b}
\end{figure}

\begin{figure}[h!]
  \centering
  \includegraphics[width=.75\linewidth]{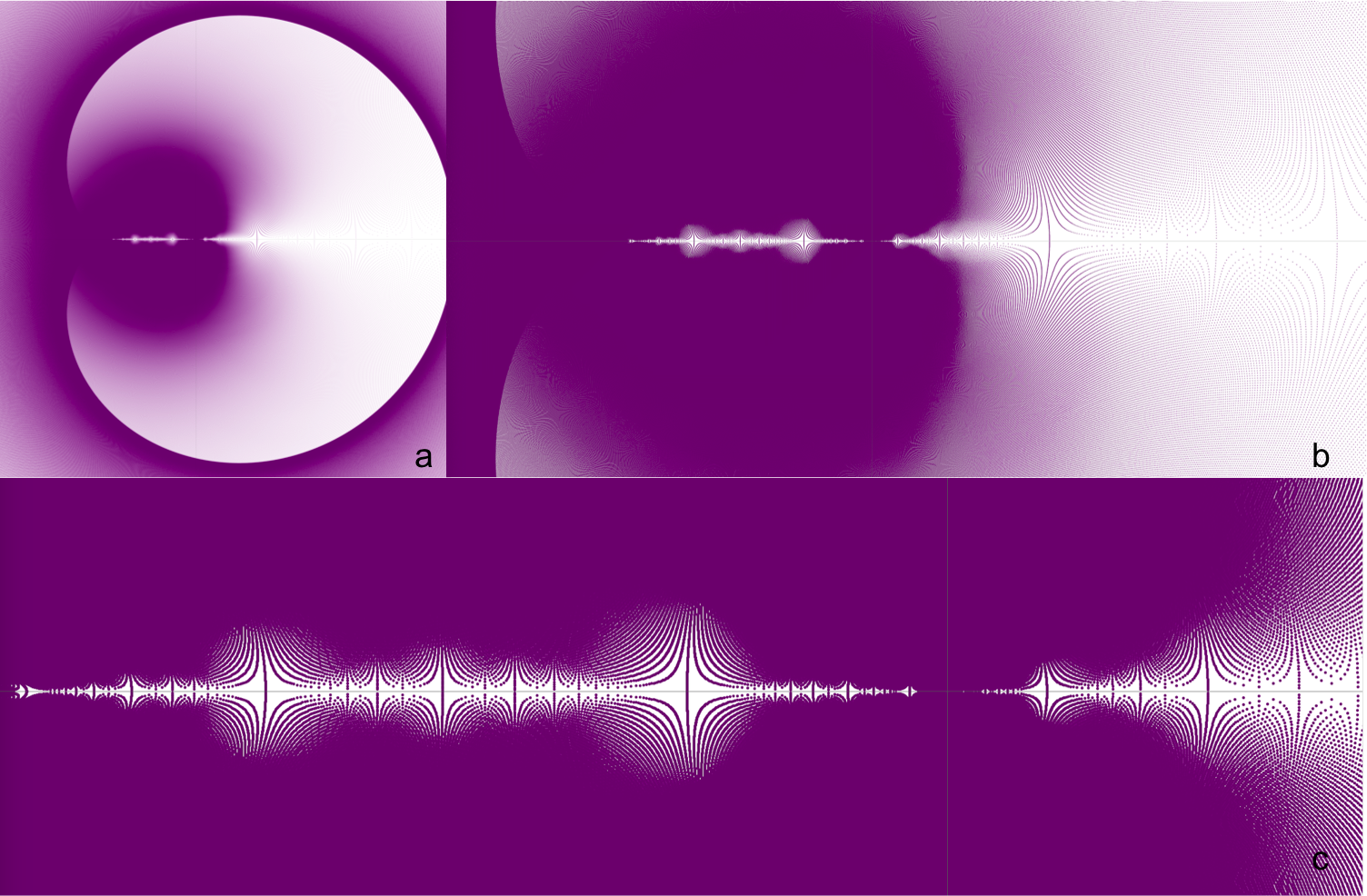}
  \caption{Appearance of the fractal pattern at 
  the 34th iteration of the Fibonacci chain.}
  \label{f2}
\end{figure}

\begin{figure}[h!]
  \centering
\includegraphics[width=0.7\linewidth]{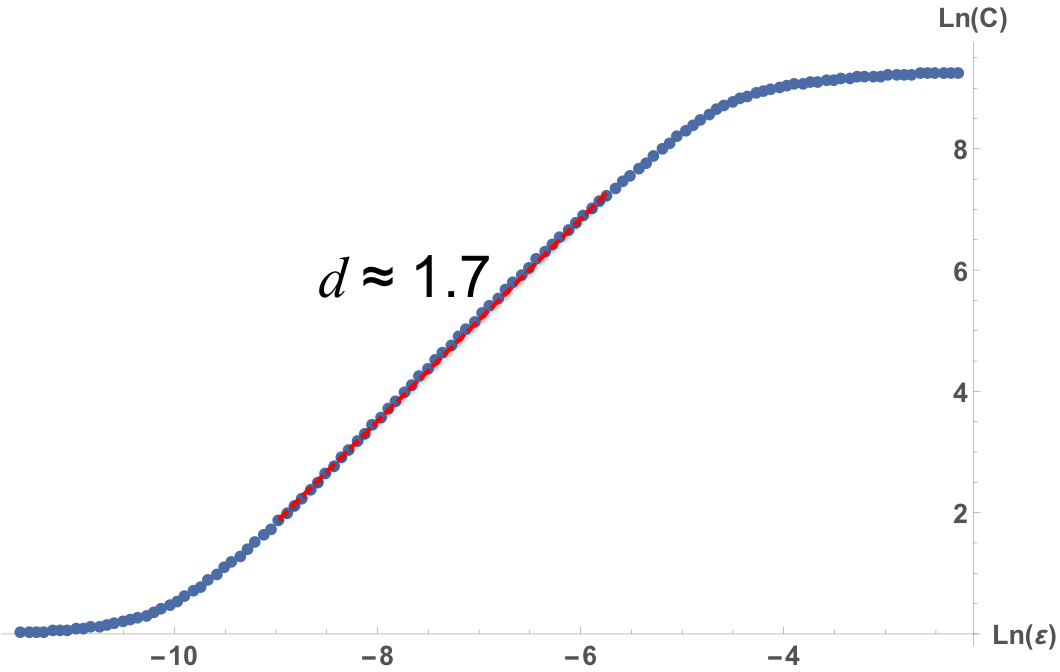}
  \caption{Pointwise dimension of the of the Fibonacci chain in Fourier space.}
  \label{f3}
\end{figure}

\section{The derivatives of the Fibonacci chain in Fourier space}

\begin{figure}[h!]
  \centering
  \includegraphics[width=0.9\textwidth,clip]{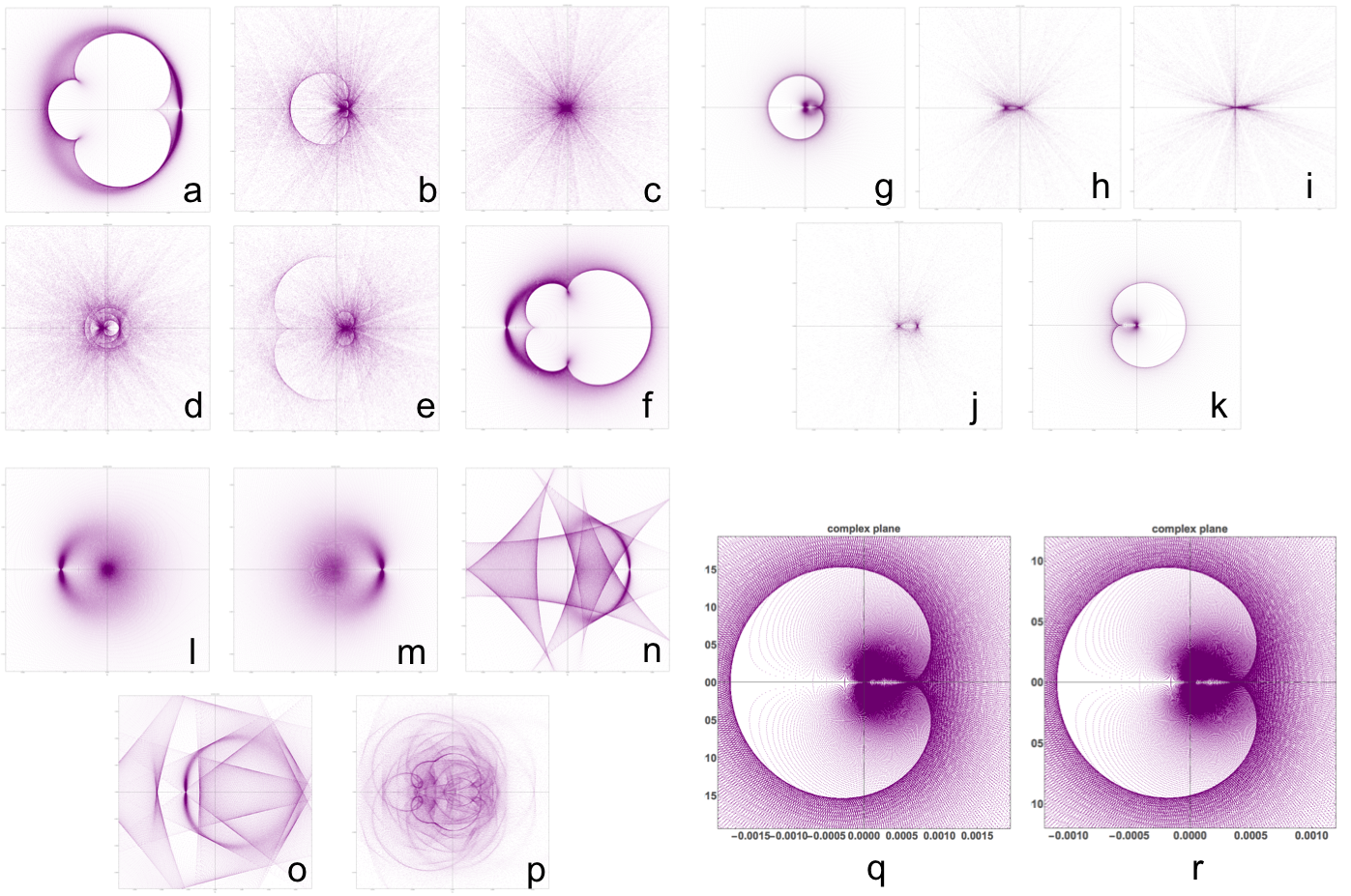}
  \caption{
The Fourier space representation of a Fibonacci chain of 26th iteration, with a) the first segment removed; b) the first two segments removed; c) the first three segments removed; d) the first six segments removed, e) the first seven segments removed, f) the first segment length replaced with 0; g) no modification; h) the last segment removed, i) the last two segments remove; j) the last 46367 segments removed; k) the last 46368 segments removed and the Fibonacci chain of 26th iteration is truncated to the Fibonacci chain of 25th iteration; l) the 1st segment $L$ replaced with $S$; m) the last two segments flipped order; n) the order of the last five segments scrambled; o) the order of the last ten segments scrambled; p) the order of the last 100 segments scrambled, and the comparison between q) the original Fibonacci chain and r) the chain with modified $L/S$ ratio where $L/S=2$.}
  \label{f4}
\end{figure}

It is tempting to conjecture that this kind of Fractal signature will appear in any finite portion of the Fibonacci chain or a Fibonacci chain that is nearly perfect (a Fibonacci chain that is generated strictly from the substitution rules mentioned earlier). A series of tests has been conducted to verify this conjecture. The results are surprising. Figure \ref{f4} shows the results of the Fourier space of the Fibonacci chain of 26th iteration with various, mostly minor modifications, such as, small truncation from the head or tail of the chain, scrambling the order of a very small part of the chain or changing of the $L/S$ ratio. These results show that the fractal pattern, especially the cardioid shape, is very sensitive to any modification except changing the $L/S$ ratio, in which case only the scaling of the fractal pattern changes. In other words, the fractal pattern is a direct result of the substitution rule (or inflation/deflation rule) which gives the Fibonacci ordering of the $L$ and $S$. The exact value of the $L$ and $S$ does not matter for the pattern. Even the smallest breakdown of the rule will result in a non-prefect-closure of the cardioid and therefore result in a completely different pattern. It is similar to the butterfly effect in chaos theory.

\begin{figure}[h!]
  \centering
  \includegraphics[width=.5\linewidth]{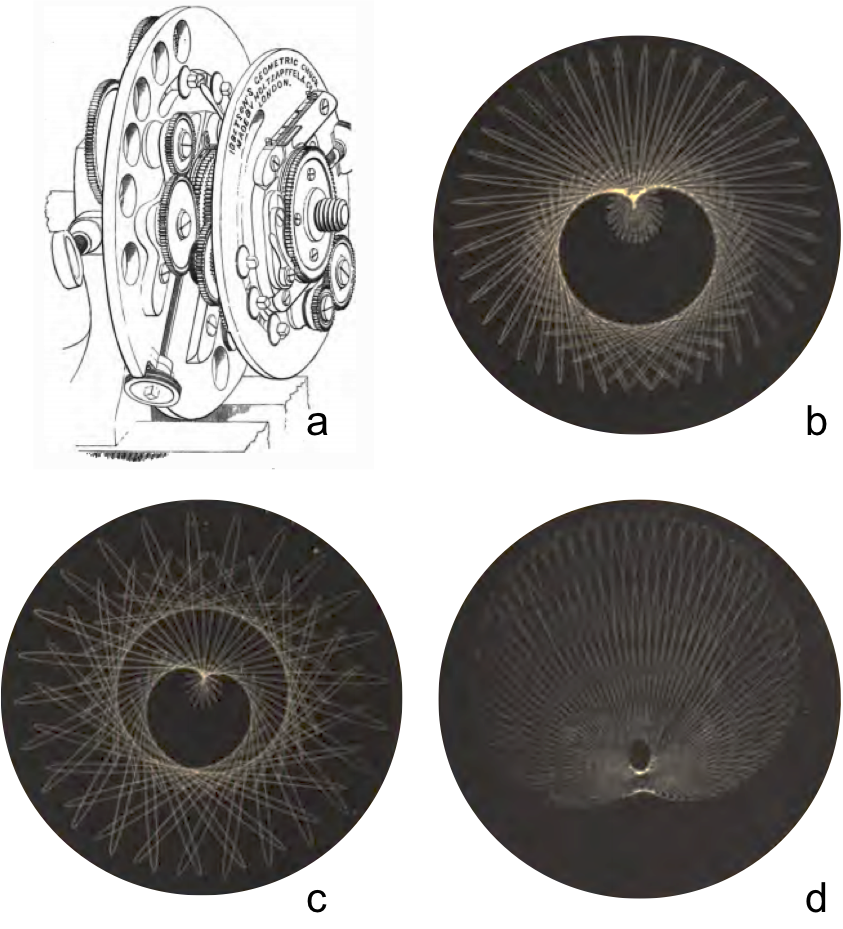}
   \caption{a) A picture of a geometric chuck; b-d) Cycloids generated by the geometric chuck with different settings.}

  \label{f5}
\end{figure}

Its worth noting that derivative patterns resemble the cycloidal pattern drawn by a geometric chuck. Figure \ref{f5} shows pictures taken from a book published 100 years ago by  Bazley \cite{Bazley1875}. Figure \ref{f5}a is a picture of a geometric chuck. Figure \ref{f5}b-d are a few cycloids created by the geometric chuck with different settings on the chuck. The fractal signature of a perfect Fibonacci chain may correspond to some perfect quasiperiodic gearing between the plates of the geometric chuck.

\section{Summary}
This paper reports a novel fractal analysis of a Fibonacci chain in Fourier space by computing the pointwise dimensions and visual inspection of the Fourier pattern. The cardioid structure in the Mandelbrot set is observed for the case of the one dimensional Fibonacci chain, with alternating mirror symmetry according to odd or even instances of the substitution iteration that generates the Fibonacci chains. The Fourier pattern is very sensitive to how strictly the chain follows the substitution rules. Any small variation in that will completely change the pattern. In comparison, the change in the $L/S$ ratio has no influence in the Fourier pattern other than a change in scaling.

\clearpage

\section{References}

\bibliographystyle{ieeetr}
\bibliography{ICQ13Fang-revised.bbl}

\begin{thebibliography}{10}

\bibitem{shechtman1984metallic}
D.~Shechtman, I.~Blech, D.~Gratias, and J.~W. Cahn, ``Metallic phase with
  long-range orientational order and no translational symmetry,'' {\em Physical
  Review Letters}, vol.~53, no.~20, p.~1951, 1984.

\bibitem{GrimmQCHamiltonian}
U.~Grimm and M.~Schreiber, ``Energy spectra and eigenstates of quasiperiodic
  tight-binding hamiltonians,'' {\em arXiv:cond-mat/0212140}, 2002.

\bibitem{RepetowiczPTHamiltonians}
P.~Repetowicz, U.~Grimm, and M.~Schreiber, ``Exact eigenstates of tight-binding
  hamiltonians on the penrose tiling,'' {\em arXiv:cond-mat/9805321}, 1998.

\bibitem{MaceFiChainFractal2016}
N.~Mac{\'e}, A.~Jagannathan, and F.~Pi{\'e}chon, ``Fractal dimensions of the
  wavefunctions and local spectral measures on the fibonacci chain,'' {\em
  arXiv:1601.00532}, 2016.

\bibitem{Ramachandrarao2000}
P.~Ramachandrarao, A.~Sina, and D.~Sanyal, ``On the fractal nature of {Penrose}
  tiling,'' {\em Current Science}, vol.~79, no.~3, pp.~364--367, 2000.

\bibitem{Yudin2001}
V.~Yudin and Y.~Karygina, ``Fractal images of quasicrystals as an example of
  {Penrose} lattice,'' {\em Crystallogrphy reports}, vol.~46, no.~6,
  pp.~922--926, 2001.

\bibitem{Fang2015}
F.~Fang and K.~Irwin, ``An icosahedral quasicrystal as a golden modification of
  the icosagrid and its connection to the e8 lattice,'' {\em arXiv:1511.07786},
  2015.

\bibitem{levine1986quasicrystals}
D.~Levine and P.~J. Steinhardt, ``Quasicrystals. i. definition and structure,''
  {\em Physical Review B}, vol.~34, no.~2, p.~596, 1986.

\bibitem{douady1984etude}
A.~Douady, J.~H. Hubbard, and P.~Lavaurs, ``Etude dynamique des polyn{\^o}mes
  complexes,'' {\em Universit{\'e} de Paris-Sud, D{\'e}p. de Math{\'e}matique},
  1984.

\bibitem{DevaneyBU}
R.~L. Devaney, http://math.bu.edu/DYSYS/FRACGEOM2, 2016.

\bibitem{Strogatz1994}
S.~H. Strogatz, {\em Nonlinear Dynamics and Chaos with applications to Physics,
  Biology, Chemistry, and Engineering}.
\newblock Perseus Books, 1994.

\bibitem{Bazley1875}
T.~S. Bazley, {\em Index to Geometric Chuck}.
\newblock Waterlow and Sons, 1875.

\end{thebibliography}

\end{document}